# Small-Gain-Based Boundary Feedback Design for Global Exponential Stabilization of 1-D Semilinear Parabolic PDEs

**Iasson Karafyllis\* and Miroslav Krstic\*\*\***

\*Dept. of Mathematics, National Technical University of Athens,
Zografou Campus, 15780, Athens, Greece, email: iasonkar@central.ntua.gr

\*\*\*Dept. of Mechanical and Aerospace Eng., University of California, San Diego, La Jolla, CA 92093-0411, U.S.A., email: krstic@ucsd.edu

**Abstract**

This paper presents a novel methodology for the design of boundary feedback stabilizers for 1-D, semilinear, parabolic PDEs. The methodology is based on the use of small-gain arguments and can be applied to parabolic PDEs with nonlinearities that satisfy a linear growth condition. The nonlinearities may contain nonlocal terms. Two different types of boundary feedback stabilizers are constructed: a linear static boundary feedback and a nonlinear dynamic boundary feedback. It is also shown that there are fundamental limitations for feedback design in the parabolic case: arbitrary gain assignment is not possible by means of boundary feedback. An example with a nonlocal nonlinear term illustrates the applicability of the proposed methodology.

**Keywords:** Parabolic PDEs, ISS, boundary feedback, feedback stabilization.

## 1. Introduction

Stabilization of 1-D parabolic Partial Differential Equations (PDEs) by means of boundary feedback control is a challenging problem that has attracted the interest of many researchers. Various methodologies for boundary feedback design in linear 1-D parabolic PDEs are available in the literature (see [1,6] and references therein), including most recently backstepping (see [17,24]). Backstepping has also been used for adaptive control of linear 1-D parabolic PDEs: see [25,16].

The motivation for the extension of the existing boundary feedback design methodologies to the case of nonlinear parabolic PDEs with nonlocal terms is strong. Nonlinear and possibly nonlocal parabolic PDEs arise in many physical problems: see for instance [2,3,15,22,23]. More specifically, a nonlinear and possibly nonlocal PDE may be an equivalent description of a system of parabolic-elliptic PDEs. Systems of parabolic-elliptic PDEs have been used in many applications: see for instance [28] for lithium-ion battery systems and [5] as well as Chapter 6 in [10] (pages 217-218) for groundwater flow. The stabilization problem of nonlinear parabolic PDEs is studied in [4], while backstepping has been extended to the case of Volterra nonlinearities in [26,27] (but see also [18]). Most of the existing stabilization results for nonlinear PDEs are local.

The purpose of the present work is the development of *global* stabilization results for the boundary feedback design problem in 1-D semilinear parabolic PDEs that may contain nonlocal terms. The methodology followed in the present work is very different from the existing boundary feedback design methodologies. The results in this paper are proved by using small-gain arguments and the Input-to-State Stability (ISS) property. The use of ISS for the study of PDEs has recently attracted the interest of the control community (see [7,11,12,13,19,20]) and use of small-gain arguments for PDEs was made in [14,15].

The focus is on the case of Dirichlet actuation at one end of the domain and Dirichlet boundary condition at the other end, but the results can be extended to other cases (Neumann actuation, Robin or Neumann boundary conditions at the non-actuated end of the domain). In order to obtain global stabilization results, we need to impose a linear growth condition on the nonlinear (and possibly nonlocal) term with restricted linear growth coefficient. Growth conditions for parabolic PDEs are also encountered in [8] for the study of global existence of solutions and lack of global controllability is proved in [9] for nonlinearities that grow faster than $|u|(\ln(1+|u|))^2$. In this work, two different boundary feedback stabilizers are provided:

- a linear static boundary feedback stabilizer, which can handle uncertain nonlinearities (Proposition 3.2), and
- a nonlinear dynamic boundary feedback stabilizer (Section 4).

Both controllers cannot handle nonlinearities which satisfy a linear growth condition with arbitrarily large linear growth coefficient. This is expected, because it is also shown that there are fundamental limitations for boundary feedback design in the parabolic case: arbitrary gain assignment is not possible (Theorem 3.3) by means of boundary feedback. The class of uncertain nonlinearities that can be handled by a linear static boundary feedback stabilizer has to satisfy a demanding linear growth condition with a small linear growth coefficient. On the other hand, this demanding linear growth condition can be avoided in certain cases and we can allow larger linear growth coefficients by using the proposed nonlinear dynamic boundary feedback stabilizer (Theorem 4.1 and Corollary 4.2). More specifically, the proposed nonlinear dynamic boundary feedback works for nonlinearities that can be expressed by linear combinations of "separable" terms of the form $\sin(\omega x)K(u)$, where $x$ is the spatial variable, $K(u)$ is a functional of the state and $\omega > 0$ is a sufficiently large constant ("wiggly in space" nonlinearity). The design of the dynamic boundary stabilizer is explicit and is based on a convenient methodology that does not require the solution of any equations. It consists of three steps:

1) We design a linear boundary feedback law that stabilizes the linear part of the PDE.
2) We design a nonlinear dynamic boundary feedback law that deals exclusively with the nonlinear and nonlocal term.
3) We combine both controllers.

It is important to notice that the first two steps are *independent* of each other. Moreover, the linear boundary feedback law may be designed by using any methodology of boundary feedback design for linear parabolic PDEs. Therefore, the proposed dynamic boundary feedback stabilizer builds on the existing design methodologies for linear PDEs and extends their applicability to the nonlinear and nonlocal case.

The structure of the present work is as follows. Section 2 presents a motivating example that cannot be handled by any other existing methodology. The example shows that very simple nonlinear and nonlocal terms may destabilize a parabolic PDE. Section 3 develops static boundary feedback stabilizers and reveals the underlying fundamental limitations for the parabolic case. Section 4 is devoted to the presentation of the proposed dynamic boundary feedback stabilizer. Section 5 provides the proofs of all results and Section 6 revisits the motivating example that was shown in Section 2. It is shown that the use of the proposed dynamic boundary feedback stabilizer can guarantee global exponential stability of the equilibrium point. The concluding remarks of the present work are given in Section 7.

Finally, it should be noted that no existence/uniqueness result for the solutions of the closed-loop system is provided in the present work. In general, the user has to assume additional regularity conditions for the nonlinear term in order to be able to guarantee existence/uniqueness of solutions for the closed-loop system using standard results (e.g., results in [21]).

**Notation.** Throughout this paper, we adopt the following notation.

* $\Re_+ := [0,+\infty)$. $Z_+$ denotes the set of non-negative integers.
* Let $A \subseteq \Re^n$ be an open set and let $\Omega \subseteq \Re$ and $A \subseteq U \subseteq \bar{A}$ be given sets. By $C^0(U)$ (or $C^0(U;\Omega)$), we denote the class of continuous mappings on $U$ (which take values in $\Omega$). By $C^k(U)$ (or $C^k(U;\Omega)$), where $k \geq 1$, we denote the class of continuous functions on $U$, which have continuous derivatives of order $k$ on $U$ (and also take values in $\Omega$).
* $L^2(0,1)$ denotes the equivalence class of measurable functions $f:[0,1] \to \Re$ for which
$$\|f\| = \left(\int_0^1 |f(x)|^2 dx\right)^{1/2} < +\infty.$$
* Let $u: \Re_+ \times [0,1] \to \Re$ be given. We use the notation $u[t]$ to denote the profile at certain $t \geq 0$, i.e., $(u[t])(x) = u(t,x)$ for all $x \in [0,1]$. $C^0\left(\Re_+; L^2(0,1)\right)$ denotes the class of continuous mappings $\Re_+ \ni t \to u[t] \in L^2(0,1)$. When $k(x)$ is differentiable with respect to $x \in [0,1]$, we use the notation $k'(x)$ for the derivative of $k$ with respect to $x \in [0,1]$, i.e., $k'(x) = \frac{dk}{dx}(x)$.

## 2. A Motivating Example

Consider the control system

$$u_t(t,x) = p u_{xx}(t,x) - q u(t,x) + (f(u[t]))(x), \text{ for } t > 0, \ x \in (0,1) \quad (2.1)$$

with

$$\begin{aligned} u(t,0) &= 0 \\ u(t,1) &= U(t) \end{aligned} \quad (2.2)$$

where $u$ is the state, $p > 0$ is the diffusion coefficient, $q \in \Re$ is the reaction coefficient, $f: L^2(0,1) \to L^2(0,1)$ is a continuous mapping with $f(0) = 0$ and $U$ is the control input.

The stability properties of the equilibrium point $0 \in L^2(0,1)$ for the open-loop system (2.1), (2.2) with $U(t) \equiv 0$ depend heavily on the values of the parameters $p, q$ and the nature of the mapping $f: L^2(0,1) \to L^2(0,1)$. For example, when

$$p = 1, \ q = 0, \ (f(u))(x) = A \sin(\omega x)\|u\| \text{ for } x \in [0,1] \quad (2.3)$$

where $A \in \Re$, $\omega > 0$ are constants then simulations indicate that the equilibrium point $0 \in L^2(0,1)$ becomes unstable for large values of $A \in \Re$. For $\omega = 20$, instability arises when $A > 447$; see also Figure 1, which shows the evolution of $\|u[t]\|$ for $\omega = 20$ and $A = 500$.

Indeed, Proposition 3.2 in the following section and results in [15] guarantee global exponential stability of the equilibrium point $0 \in L^2(0,1)$ for the open-loop system (2.1), (2.2), (2.3) with $U(t) \equiv 0$, only when

$$|A| < \frac{2\pi^2 \sqrt{\omega}}{\sqrt{2\omega - \sin(2\omega)}} \quad (2.4)$$

It should be noted here that (2.4) is a conservative estimation of the stability region: for $\omega = 20$ inequality (2.4) requires $|A| < 14.08957$, while simulations indicate that the equilibrium point $0 \in L^2(0,1)$ is stable for $|A| \leq 430$.

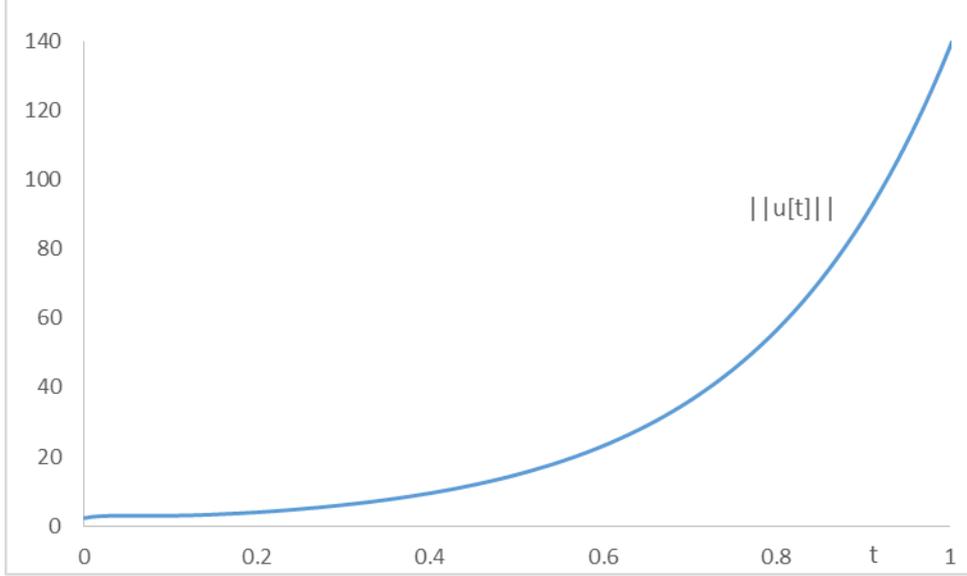

**Fig.1:** The evolution of $\|u[t]\|$ for the open-loop system (2.1), (2.2), (2.3) with $U(t) \equiv 0$, $\omega = 20$ and $A = 500$. The initial condition is

$$u_0(x) = \sqrt{2}\sin(\pi x) + 2\sqrt{2}\sin(2\pi x) + \frac{\sqrt{2}}{10}\sin(3\pi x) + \frac{\sqrt{2}}{100}\sum_{n=37}^{42}\sin(n\pi x), \quad x \in [0,1].$$

It is clear that the stabilization of the equilibrium point $0 \in L^2(0,1)$ by means of a boundary feedback law is required for large values of $A \in \Re$. However, there are no available results in the literature that can handle the nonlinear and nonlocal term $(f(u))(x) = A\sin(\omega x)\|u\|$. This term causes excitement of all modes of the state $u[t]$ when $\frac{\omega}{\pi} \notin Z_+$ and this creates an exponential increase of the norm of the state $\|u[t]\|$ (see Fig.1).

## 3. Small-Gain-Based Static Boundary Feedback Design

Consider the control system

$$u_t(t,x) = pu_{xx}(t,x) - qu(t,x) + v(t,x), \text{ for } t > 0, \; x \in (0,1) \quad (3.1)$$

with

$$\begin{aligned} u(t,0) &= 0 \\ u(t,1) &= U(t) \end{aligned}, \text{ for } t > 0 \quad (3.2)$$

where $u$ is the state, $p > 0$ is the diffusion coefficient, $q \in \Re$ is the reaction coefficient, $v : \Re_+ \times [0,1] \to \Re$ is a distributed input and $U$ is the control input. The distributed input $v$ is going to be used later for the quantification of the effect of nonlinear and possibly non-local terms that may appear in the right-hand side of a semilinear, 1-D parabolic PDE.

**Definition 3.1 (Input-to-State Stabilization):** *We say that the kernel $k \in C^0([0,1])$ achieves Input-to-State Stabilization with gain $\frac{\gamma}{p\pi^2} \geq 0$ of system (3.1), (3.2), if there exist constants $G, \sigma > 0$ so that every solution $u \in C^0(\Re_+; L^2(0,1)) \cap C^1((0,+\infty) \times [0,1])$ with $u[t] \in C^2([0,1])$ for $t > 0$ of the closed-loop system (3.1), (3.2) with*

$$U(t) = \int_0^1 k(x) u(t,x) dx, \text{ for } t > 0 \tag{3.3}$$

*satisfies the estimate:*

$$\|u[t]\| \leq G \exp(-\sigma t) \|u[0]\| + \frac{\gamma}{p\pi^2} \sup_{0 \leq s \leq t} (\|v([s])\|), \text{ for all } t \geq 0 \tag{3.4}$$

The motivation for the study of kernels that achieve Input-to-State Stabilization of system (3.1), (3.2) comes from the following result.

**Proposition 3.2 (Robust Global Exponential Stabilization of an Uncertain Nonlinear Parabolic PDE):** *Consider the control system*

$$u_t(t,x) = p u_{xx}(t,x) - q u(t,x) + d(t)(f(u[t]))(x), \text{ for } t > 0, \ x \in (0,1) \tag{3.5}$$

*with (3.2), where $u$ is the state, $p > 0$ is the diffusion coefficient, $q \in \Re$ is the reaction coefficient, $f : L^2(0,1) \to L^2(0,1)$ is a continuous mapping, $d : \Re_+ \to [-1,1]$ is a disturbance (vanishing perturbation) and $U$ is the control input. Suppose that the kernel $k \in C^0([0,1])$ achieves Input-to-State Stabilization with gain $\frac{\gamma}{p\pi^2} > 0$ of system (3.1), (3.2). Moreover, suppose that there exists a constant $M \geq 0$ such that the following linear growth condition holds:*

$$\|f(u)\| \leq M \|u\|, \text{ for all } u \in L^2(0,1) \tag{3.6}$$

*Finally, suppose that the following small-gain condition holds for $M$:*

$$M < \gamma^{-1} p\pi^2 \tag{3.7}$$

*Then the kernel $k \in C^0([0,1])$ achieves robust global exponential stabilization of system (3.5), (3.2), i.e., there exist constants $R, \omega > 0$ so that every solution $u \in C^0(\Re_+; L^2(0,1)) \cap C^1((0,+\infty) \times [0,1])$ with $u[t] \in C^2([0,1])$ for $t > 0$ of the closed-loop system (3.5), (3.2) with (3.3) satisfies the estimate:*

$$\|u[t]\| \leq R \exp(-\omega t) \|u[0]\|, \text{ for all } t \geq 0 \tag{3.8}$$

Unfortunately, we cannot achieve robust global exponential stabilization of system (3.5), (3.2) for nonlinearities with arbitrarily large coefficient $M$ of linear growth. The reason is the following result, which poses a fundamental limitation on the achievable gain.

**Theorem 3.3 (Gain Assignment is not Possible by Means of Boundary Feedback):** *If the kernel $k \in C^0([0,1])$ achieves Input-to-State Stabilization with gain $\frac{\gamma}{p\pi^2 + q} \geq 0$ of system (3.1), (3.2) with $q \geq 0$ then*

$$\gamma \geq \sup_{m=1,2,\ldots} \left(\sqrt{g_m(\mu)}\right) \tag{3.9}$$

where $\mu = \frac{1}{\pi}\sqrt{\frac{q}{p}}$ and

$$g_m(0) := \max\left\{\sum_{n=1}^{m}\frac{v_n^2}{n^4} - \frac{6}{\pi^2}\left(\sum_{n=1}^{m}\frac{(-1)^n v_n}{n^3}\right)^2 : (v_1,\ldots,v_m) \in \Re^m, \sum_{n=1}^{m} v_n^2 = 1\right\} \quad (3.10)$$

$g_m(\mu) :=$

$$\max\left\{\sum_{n=1}^{m} v_n^2\left(\frac{1+\mu^2}{n^2+\mu^2}\right)^2 - \frac{8}{\pi^2}\left(\frac{\mu\pi\sinh^2(\mu\pi)}{\sinh(2\mu\pi)-2\mu\pi}\right)\left(\sum_{n=1}^{m}\frac{(-1)^n(1+\mu^2)nv_n}{(n^2+\mu^2)^2}\right)^2 : (v_1,\ldots,v_m) \in \Re^m, \sum_{n=1}^{m} v_n^2 = 1\right\},$$

for $\mu > 0$ (3.11)

**Remark:** Notice that for every $\mu \geq 0$ the sequence $g_m(\mu)$ is non-decreasing and bounded. Therefore, the least upper bound appearing in (3.9) is well-defined. Moreover, by direct computation of the right hand sides of (3.10), (3.11) we have:

$$g_1(0) = 1 - \frac{6}{\pi^2}, \quad g_2(0) = \frac{3\sqrt{(10\pi^2-63)^2+256}+34\pi^2-195}{64\pi^2}$$

$$g_1(\mu) = 1 - \frac{8\mu\sinh^2(\mu\pi)}{(1+\mu^2)^2\pi(\sinh(2\mu\pi)-2\mu\pi)}, \text{ for } \mu > 0$$

Notice that the zero kernel achieves Input-to-State Stabilization with gain $\frac{1}{p\pi^2+q}$ of system (3.1), (3.2) (see [15]). The following result provides an upper bound for the ISS gain that is achieved by the linear single-mode boundary feedback $U(t) = -\pi r\int_0^1 u(t,x)\sin(\pi x)dx$ with $r \geq 0$.

**Theorem 3.4 (Achievable ISS Gain by Linear Single-Mode Boundary Feedback):** *Consider the control system (3.1), (3.2) with $q \geq 0$. Then, for every $r \geq 0$, every solution $u \in C^0(\Re_+; L^2(0,1)) \cap C^1((0,+\infty) \times [0,1])$ with $u[t] \in C^2([0,1])$ for $t > 0$ of the closed-loop system (3.1), (3.2) with*

$$U(t) = -\pi r \int_0^1 u(t,x)\sin(\pi x)dx, \text{ for } t > 0 \quad (3.12)$$

*satisfies the following ISS estimate*

$$\|u[t]\| \leq \exp\left(-(p\pi^2+q)t\right)\sqrt{2(1+\varepsilon^{-1})}\left(1+r\sqrt{\frac{8(\pi^2-6)}{27}}\right)\|u[0]\|$$

$$+ \sqrt{\frac{(1+\varepsilon)(\mu^2+\omega^*)}{(\mu^2+\omega^*-\sigma)}}\frac{b(r,\mu)}{p\pi^2+q}\sup_{0\leq s\leq t}\left(\exp(-\sigma p\pi^2(t-s))\|v[s]\|\right),$$

for all $t \geq 0$, $\varepsilon > 0$, $\sigma \in \left[0, \omega^* + \mu^2\right)$ (3.13)

*where $\mu := \frac{1}{\pi}\sqrt{\frac{q}{p}}$,*

$$b(r,\mu) := \frac{1+\mu^2}{2}\sqrt{\min\left\{\frac{L(r,\omega,\lambda)}{(\mu^2+\omega)(4-\omega)} : \lambda \in (0,1), \omega \in \left(-\mu^2, \min(4, 1+r)\right)\right\}}, \quad (3.14)$$

$$L(r,\omega,\lambda) := \lambda^{-1}\max\left(\frac{\lambda(4-\omega)}{1+r-\omega}\left(1+r^2\sum_{n=2}^{\infty}\frac{n^2}{(n^2-\omega)\left(n^2-\omega-\lambda(4-\omega)\right)}\right), 1\right), \quad (3.15)$$

and $\omega^* \in \left(-\mu^2, \min(4, 1+r)\right)$, $\lambda^* \in (0,1)$ are constants that satisfy $\left(\frac{2b(r,\mu)}{1+\mu^2}\right)^2 = \frac{L(r,\omega^*,\lambda^*)}{(\mu^2+\omega^*)(4-\omega^*)}$.

Notice that since $\varepsilon > 0$ is arbitrary, it follows from (3.13) that the achievable gain by the linear single-mode boundary feedback (3.12) with $r \geq 0$ is any number greater than $\frac{b(r,\mu)}{p\pi^2 + q}$. For $\mu = 0$ (i.e., $q = 0$) we may use the inequality $L(r,\omega,\lambda) \leq \tilde{L}(r,\omega,\lambda) := \lambda^{-1}\max\left(\frac{\lambda(4-\omega)}{1+r-\omega}\left(1+r^2\frac{8(\pi^2-6)}{3(4-\omega)^2(1-\lambda)}\right), 1\right)$ in order to approximate $b(r,0)$ from above. Indeed, we have $b(r,0) \leq ub(r) := \frac{1}{2}\sqrt{\min\left\{\frac{\tilde{L}(r,\omega,\lambda)}{\omega(4-\omega)} : \lambda \in (0,1), \omega \in \left(0, \min(4, 1+r)\right)\right\}}$. Figure 2 shows the graph of the function $ub(r)$ (the upper bound of $b(r,0)$). The smallest gain $\frac{ub(r)}{p\pi^2}$ is achieved for $r = 0.91$ with $ub(r)$ being equal to $0.699268$. Indeed, $g_2(0)$ as defined by (3.10) is equal to $0.639003$ and consequently, the linear single-mode boundary feedback $U(t) = -0.91\pi \int_0^1 u(t,x)\sin(\pi x)dx$ is close to the lowest possible gain predicted by (3.9).

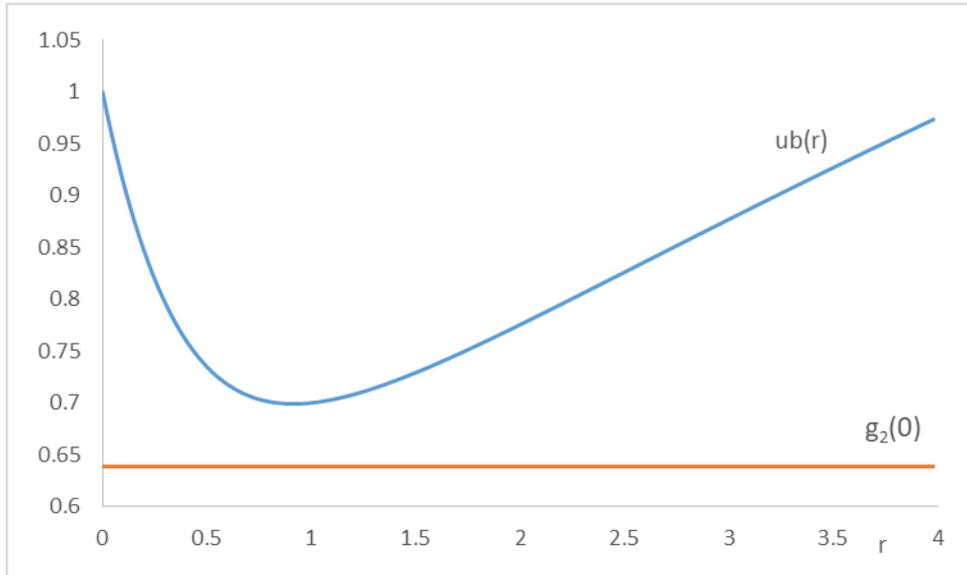

**Fig. 2:** The achievable gain for $q=0$ by the linear single-mode boundary feedback
$$U(t) = -\pi r \int_0^1 u(t,x)\sin(\pi x)dx \text{ with } r \geq 0.$$

# 4. Small-Gain-Based Dynamic Boundary Feedback Design

Can we handle nonlinearities with linear-growth coefficients that violate the small-gain condition (3.7)? The answer is "yes" but in this case we have to eliminate the uncertainty. This is guaranteed by the following result, which utilizes a dynamic nonlinear boundary feedback law.

**Theorem 4.1 (Exponential Stabilization by Means of Dynamic Nonlinear Boundary Feedback):** *Suppose that the kernel $k \in C^0([0,1])$ achieves Input-to-State Stabilization with gain $\frac{\gamma}{p\pi^2} \geq 0$ of system (3.1), (3.2). Consider the control system*

$$u_t(t,x) = p u_{xx}(t,x) - q u(t,x) + (f(u[t]))(x), \text{ for } t > 0, \ x \in (0,1) \tag{4.1}$$

*with (3.2), where $u$ is the state, $p > 0$ is the diffusion coefficient, $q \in \Re$ is the reaction coefficient, $f : L^2(0,1) \to L^2(0,1)$ is a continuous mapping and $U$ is the control input. Suppose that there exist functions $\varphi_i \in C^2([0,1])$ with $\varphi_i(0) = 0$ ($i = 1,...,n$), continuous functionals $K_i : L^2(0,1) \to \Re$ ($i = 1,...,n$), constants $M \geq 0$, $P_i \geq 0$, $\omega_i > 0$ ($i = 1,...,n$), so that the following growth conditions hold:*

$$|K_i(u)| \leq P_i \|u\|, \text{ for } i = 1,...,n \text{ and } u \in L^2(0,1) \tag{4.2}$$

$$\left\| f(u) - \sum_{i=1}^{n} \varphi_i K_i(u) \right\| \leq M \|u\|, \text{ for all } u \in L^2(0,1) \tag{4.3}$$

*Furthermore suppose that the following small-gain condition is satisfied by all $P_i \geq 0$ ($i = 1,...,n$):*

$$\frac{\gamma}{p\pi^2} \left( M + \sum_{i=1}^{n} \omega_i^{-2} P_i \left\| p \varphi_i'' - (q - \omega_i^2) \varphi_i \right\| \right) + \sum_{i=1}^{n} \omega_i^{-2} P_i \|\varphi_i\| < 1 \tag{4.4}$$

*Then $0 \in L^2(0,1) \times \Re^n$ is globally exponentially stable for the closed-loop system (4.1), (3.2) with*

$$U(t) = \sum_{i=1}^{n} \varphi_i(1) \xi_i(t) + \int_0^1 k(x) u(t,x) dx - \sum_{i=1}^{n} \xi_i(t) \int_0^1 \varphi_i(x) k(x) dx \tag{4.5}$$

$$\dot{\xi}_i(t) = -\omega_i^2 \xi_i(t) + K_i(u[t]), \ i = 1,...,n \tag{4.6}$$

*i.e., there exist constants $R, \delta > 0$ so that every solution $u \in C^0(\Re_+; L^2(0,1)) \cap C^1((0,+\infty) \times [0,1])$, $\xi_i \in C^0(\Re_+) \cap C^1((0,+\infty))$, $i = 1,...,n$ with $u[t] \in C^2([0,1])$ for $t > 0$ of the closed-loop system (4.1), (3.2) with (4.5), (4.6) satisfies the estimate:*

$$\|u[t]\| + \sum_{i=1}^{n} |\xi_i(t)| \leq R \exp(-\delta t) \left( \|u[0]\| + \sum_{i=1}^{n} |\xi_i(0)| \right), \text{ for all } t \geq 0 \tag{4.7}$$

**Remarks: (a)** When $q \geq 0$ and $k(x) = -\pi r \sin(\pi x)$ with $r \geq 0$ (single-mode boundary feedback), then Theorem 3.4 guarantees that the small-gain condition (4.4) takes the form

$$\frac{b(r,\mu)}{p\pi^2 + q} \left( M + \sum_{i=1}^{n} \omega_i^{-2} P_i \left\| p \varphi_i'' - (q - \omega_i^2) \varphi_i \right\| \right) + \sum_{i=1}^{n} \omega_i^{-2} P_i \|\varphi_i\| < 1 \tag{4.8}$$

where $\mu := \frac{1}{\pi} \sqrt{\frac{q}{p}}$ and $b(r,\mu)$ is the function defined by (3.14).

**(b)** It should be noticed that the controller (4.5), (4.6) is a *nonlinear* dynamic boundary feedback law. Thus, the response of the closed-loop system (4.1), (3.2) with (4.5), (4.6) does not depend only on the initial condition of the distributed state $u[0]$ but depends also on the initial conditions of the internal controller states $\xi_i(0)$, $i=1,...,n$. The small-gain condition (4.4) reveals that dynamic boundary feedback stabilization by the integral controller (4.5), (4.6) is possible only for nonlinearities $f(u)$ that can be approximated closely (in the $L^2$ norm) by linear combinations of "separable" nonlinear terms, i.e., by $\sum_{i=1}^{n}\varphi_i K_i(u)$. More specifically, the nonlinear term may satisfy the following assumption.

**(H)** *There exist a countable family of continuous functionals $K_i : L^2(0,1) \to \Re$ $i=1,2,...$ and constants $P_i \geq 0$, $\omega_i > 0$, $i=1,2,...$, satisfying (4.2) and*

$$\sum_{i=1}^{\infty} \omega_i^{-2} P_i \|\varphi_i\| < 1 \tag{4.9}$$

*where the functions $\varphi_i$ satisfy $\varphi_i(x) = x$ when $q = \omega_i^2$, $\varphi_i(x) = \sin\left(x\sqrt{p^{-1}(\omega_i^2 - q)}\right)$ when $q < \omega_i^2$ and $\varphi_i(x) = \sinh\left(x\sqrt{p^{-1}(q - \omega_i^2)}\right)$ when $q > \omega_i^2$, such that for every $\varepsilon > 0$, there exists an integer $n \geq 1$ satisfying*

$$\left\| f(u) - \sum_{i=1}^{n} \varphi_i K_i(u) \right\| \leq \varepsilon \|u\|, \text{ for all } u \in L^2(0,1) \tag{4.10}$$

The fact that nonlinear (and possibly nonlocal) continuous mappings $f : L^2(0,1) \to L^2(0,1)$ which satisfy assumption (H) can be handled by Theorem 4.1 becomes apparent by the following corollary.

**Corollary 4.2 (Exponential Stabilization by Means of Dynamic Nonlinear Boundary Feedback):** *Suppose that the kernel $k \in C^0([0,1])$ achieves Input-to-State Stabilization with gain $\frac{\gamma}{p\pi^2} \geq 0$ of system (3.1), (3.2). Consider the control system (4.1) with (3.2), where $u$ is the state, $p > 0$ is the diffusion coefficient, $q \in \Re$ is the reaction coefficient, $f : L^2(0,1) \to L^2(0,1)$ is a continuous mapping and $U$ is the control input. Suppose that assumption (H) holds for the mapping $f : L^2(0,1) \to L^2(0,1)$. Then there exists an integer $n \geq 1$, functions $\varphi_i \in C^2([0,1])$ with $\varphi_i(0) = 0$ ($i=1,...,n$), continuous functionals $K_i : L^2(0,1) \to \Re$ ($i=1,...,n$), constants $\omega_i > 0$ ($i=1,...,n$), so that $0 \in L^2(0,1) \times \Re^n$ is globally exponentially stable for the closed-loop system (4.1), (3.2) with (4.5), (4.6).*

**Proof:** Simply use assumption (H) with $\varepsilon > 0$ sufficiently small so that $\frac{\gamma \varepsilon}{p\pi^2} + L < 1$, where $L := \sum_{i=1}^{\infty} \omega_i^{-2} P_i \|\varphi_i\| < 1$. The rest of the proof is a consequence of Theorem 4.1. ◁

It is clear that when assumption (H) holds, then the design of the dynamic boundary stabilizer (4.5), (4.6) is explicit and is based on a convenient methodology that does not require the solution of any equations. The methodology consists of three steps:

1) We design the kernel $k \in C^0([0,1])$ that achieves Input-to-State Stabilization of the linear part of the PDE (4.1). It should be noted that since $q \in \Re$ does not necessarily satisfy the inequality $q > -p\pi^2$, the linear part of the PDE (4.1) may be unstable.
2) We select the countable family of continuous functionals $K_i : L^2(0,1) \to \Re$ $i=1,2,...$ and constants $P_i \geq 0$, $\omega_i > 0$, $i = 1,2,...$, satisfying (4.2), (4.9), (4.10). Notice that the selection depends exclusively on the nonlinear and nonlocal term $f : L^2(0,1) \to L^2(0,1)$.
3) We combine both controllers and construct the dynamic boundary stabilizer (4.5), (4.6).

It is important to notice that the first two steps are *independent* of each other. Moreover, the kernel $k \in C^0([0,1])$ may be designed by using any available boundary feedback design methodology for linear parabolic PDEs. Particularly, when $q < 0$ then we may select a boundary feedback designed by means of backstepping (see [17,24]) and take $k(x) = -p^{-1}(r-q) x \dfrac{I_1\left(\sqrt{p^{-1}(r-q)(1-x^2)}\right)}{\sqrt{p^{-1}(r-q)(1-x^2)}}$ with $r \geq 0$, where $I_1$ is the modified Bessel function of the first kind of order one. The fact that this kernel achieves Input-to-State Stabilization of system (3.1), (3.2) follows from an analysis similar to that given in [12,15]. Therefore the proposed dynamic boundary feedback stabilizer builds on the existing design methodologies for linear PDEs and extends their applicability to the nonlinear and nonlocal case.

## 5. Proofs of Main Results

**Proof of Theorem 3.3:** The proof relies on the following fact.

**Fact I:** *Consider system (3.1), (3.2) with (3.3). If the kernel $k \in C^0([0,1])$ achieves Input-to-State Stabilization with gain $\dfrac{\gamma}{p\pi^2 + q} \geq 0$ of system (3.1), (3.2) then for every $\bar{v} \in C^1([0,1])$, for which the boundary value problem*

$$pw''(x) - qw(x) + \bar{v}(x) = 0, \text{ for } x \in [0,1] \qquad (5.1)$$

$$w(0) = w(1) - \int_0^1 k(z) w(z) dz = 0 \qquad (5.2)$$

*has a solution $w \in C^2([0,1])$, the following inequality holds*

$$\|w\| \leq \dfrac{\gamma}{p\pi^2 + q} \|\bar{v}\|. \qquad (5.3)$$

Indeed, Fact I is a direct consequence of estimate (3.4) and the fact that $u[t] \equiv w$ is a solution of system (3.1), (3.2) with (3.3), $v[t] \equiv \bar{v}$. Therefore, estimate (3.4) implies that

$$\|w\| \leq G\exp(-\sigma t)\|w\| + \frac{\gamma}{p\pi^2}\|\bar{v}\|, \text{ for all } t \geq 0$$

from which we obtain estimate (5.3).

Let $m \geq 1$ be an arbitrary integer and let arbitrary $\bar{v}_n$, $n=1,...,m$ be given. It is easy to see that for the selection

$$\bar{v}(x) = \sqrt{2}\sum_{n=1}^{m} \bar{v}_n \sin(n\pi x), \text{ for } x \in [0,1] \tag{5.4}$$

any solution of the boundary-value problem (5.1), (5.2) must be of the form:

$$w(x) = Ax + \sqrt{2}\sum_{n=1}^{m} \frac{\bar{v}_n}{pn^2\pi^2}\sin(n\pi x), \text{ for } x \in [0,1] \text{ and } q = 0 \tag{5.5}$$

$$w(x) = A\sinh(\mu\pi x) + \sqrt{2}\sum_{n=1}^{m} \frac{\bar{v}_n}{pn^2\pi^2 + q}\sin(n\pi x), \text{ for } x \in [0,1] \text{ and } q > 0 \tag{5.6}$$

where $\mu = \frac{1}{\pi}\sqrt{\frac{q}{p}}$ and $A \in \Re$ is a constant that satisfies

$$A\left(1 - \int_0^1 xk(x)dx\right) = \sqrt{2}\sum_{n=1}^{m} \frac{\bar{v}_n}{pn^2\pi^2}\int_0^1 k(x)\sin(n\pi x)dx, \text{ for } q = 0 \tag{5.7}$$

$$A\left(\sinh(\mu\pi) - \int_0^1 k(x)\sinh(\mu\pi x)dx\right) = \sqrt{2}\sum_{n=1}^{m} \frac{\bar{v}_n}{pn^2\pi^2 + q}\int_0^1 k(x)\sin(n\pi x)dx, \text{ for } q > 0 \tag{5.8}$$

The fact that both equations (5.7), (5.8) are solvable with respect to $A$ is a direct consequence of the following fact.

**Fact II:** *Consider system (3.1), (3.2) with (3.3). If the kernel $k \in C^0([0,1])$ achieves Input-to-State Stabilization with gain $\frac{\gamma}{p\pi^2 + q} \geq 0$ of system (3.1), (3.2) then the boundary value problem (5.1), (5.2) with $\bar{v} \equiv 0$ has the unique solution $w \equiv 0$.*

Fact II follows from the observation that for the disturbance-free control system (3.1), (3.2), with (3.3) and $v[t] \equiv 0$, ISS implies global exponential stability and therefore uniqueness of the equilibrium point. Any non-zero solution of the boundary-value problem (5.1), (5.2) with $\bar{v} \equiv 0$ would imply a non-zero equilibrium point for the disturbance-free control system (3.1), (3.2) with (3.3) and $v[t] \equiv 0$.

It follows from Fact II that when $q = 0$ we must have $\left(1 - \int_0^1 xk(x)dx\right) \neq 0$ (since otherwise the boundary-value problem (5.1), (5.2) with $\bar{v} \equiv 0$ would have an infinite number of solutions of the form $w(x) = Ax$ with arbitrary $A \in \Re$). Similarly, it follows from Fact II that when $q > 0$ we must have $\left(\sinh(\mu\pi) - \int_0^1 k(x)\sinh(\mu\pi x)dx\right) \neq 0$ (since otherwise the boundary-value problem (5.1), (5.2) with $\bar{v} \equiv 0$ would have an infinite number of solutions of the form $w(x) = A\sinh(\mu\pi x)$ with arbitrary $A \in \Re$). Consequently, both equations (5.7), (5.8) are uniquely solvable with respect to $A$ for every possible selection of $\bar{v}_n$, $n = 1, \ldots, m$. Therefore, for every possible selection of $\bar{v}_n$, $n = 1, \ldots, m$, the boundary-value problem (5.1), (5.2) with $\bar{v}$ given by (5.4) has a unique solution $w \in C^2([0,1])$, which is given by (5.5), (5.6), (5.7) and (5.8).

Next we rely on the following fact.

**Fact III:** *Suppose that* $w = f + Ah$, *where* $f, h \in C^2([0,1])$, $h \neq 0$ *and* $A \in \Re$. *Then*
$$\|w\|^2 \geq \|f\|^2 - \frac{\langle f, h \rangle^2}{\|h\|^2}.$$

Fact III is a direct consequence of the fact that the polynomial $g(A) := \|f\|^2 + 2A\langle f, h \rangle + A^2 \|h\|^2$ satisfies $g(A) = \|w\|^2$ and $g(A) \geq \|f\|^2 - \frac{\langle f, h \rangle^2}{\|h\|^2}$ for all $A \in \Re$.

It follows from Fact III and (5.5), (5.6) that for every possible selection of $\bar{v}_n$, $n = 1, \ldots, m$, the unique solution $w \in C^2([0,1])$ of the boundary-value problem (5.1), (5.2) with $\bar{v}$ given by (5.4) must satisfy the following inequalities (recall $\mu = \frac{1}{\pi}\sqrt{\frac{q}{p}}$):

$$\|w\|^2 \geq \frac{1}{p^2\pi^4}\sum_{n=1}^m \frac{\bar{v}_n^2}{n^4} - \frac{6}{p^2\pi^6}\left(\sum_{n=1}^m \frac{(-1)^n \bar{v}_n}{n^3}\right)^2, \text{ for } q = 0 \tag{5.9}$$

$$\|w\|^2 \geq \sum_{n=1}^m \frac{\bar{v}_n^2}{(pn^2\pi^2 + q)^2} - \frac{8\mu\pi \sinh^2(\mu\pi)}{\sinh(2\mu\pi) - 2\mu\pi}\left(\sum_{n=1}^m \frac{(-1)^n n\pi p\bar{v}_n}{(pn^2\pi^2 + q)^2}\right)^2, \text{ for } q > 0 \tag{5.10}$$

Inequality (3.9) is a direct consequence of Fact I, (5.9), (5.10), definition $\mu = \frac{1}{\pi}\sqrt{\frac{q}{p}}$, the fact that the function $\bar{v}$ given by (5.4) satisfies $\|\bar{v}\|^2 = \sum_{n=1}^m \bar{v}_n^2$ and the fact that the integer $m \geq 1$ and $\bar{v}_n$, $n = 1, \ldots, m$ are arbitrary. The proof is complete. ◁

**Proof of Proposition 3.2:** By virtue of inequality (3.7), there exists sufficiently small $\varepsilon > 0$ such that

$$\frac{\gamma(1+\varepsilon)M}{p\pi^2}<1 \tag{5.11}$$

Lemma 7.1 in [15] in conjunction with inequality (3.4) implies that there exists a constant $\delta \in (0,\sigma)$ with the following property: every solution $u \in C^0(\Re_+; L^2(0,1)) \cap C^1((0,+\infty) \times [0,1])$ with $u[t] \in C^2([0,1])$ for $t>0$ of the closed-loop system (3.1), (3.2) with (3.3) satisfies the estimate:

$$\|u[t]\| \leq G\exp(-\delta t)\|u[0]\| + \frac{\gamma(1+\varepsilon)}{p\pi^2} \sup_{0 \leq s \leq t}\left(\|v([s])\|\exp(-\delta(t-s))\right), \text{ for all } t \geq 0 \tag{5.12}$$

Notice that every solution $u \in C^0(\Re_+; L^2(0,1)) \cap C^1((0,+\infty) \times [0,1])$ with $u[t] \in C^2([0,1])$ for $t>0$ of the closed-loop system (3.5), (3.2) with (3.3) is a solution $u \in C^0(\Re_+; L^2(0,1)) \cap C^1((0,+\infty) \times [0,1])$ with $u[t] \in C^2([0,1])$ for $t>0$ of the closed-loop system (3.1), (3.2) with (3.3) and

$$v(t,x) = d(t)(f(u[t]))(x), \text{ for } t \geq 0, \ x \in [0,1] \tag{5.13}$$

or equivalently

$$v[t] = d(t)f(u[t]), \text{ for } t \geq 0 \tag{5.14}$$

Consequently, it follows from (3.6), (5.12), (5.14) and the fact that $|d(t)| \leq 1$ for all $t \geq 0$, that every solution $u \in C^0(\Re_+; L^2(0,1)) \cap C^1((0,+\infty) \times [0,1])$ with $u[t] \in C^2([0,1])$ for $t>0$ of the closed-loop system (3.5), (3.2) with (3.3) satisfies the estimate:

$$\|u[t]\|\exp(\delta t) \leq G\|u[0]\| + \frac{\gamma(1+\varepsilon)M}{p\pi^2} \max_{0 \leq s \leq t}\left(\|u[s]\|\exp(\delta s)\right), \text{ for all } t \geq 0 \tag{5.15}$$

Inequality (5.15) in conjunction with (5.11) implies that every solution $u \in C^0(\Re_+; L^2(0,1)) \cap C^1((0,+\infty) \times [0,1])$ with $u[t] \in C^2([0,1])$ for $t>0$ of the closed-loop system (3.5), (3.2) with (3.3) satisfies the estimate:

$$\max_{0 \leq s \leq t}\left(\|u[s]\|\exp(\delta s)\right) \leq G\left(1 - \frac{\gamma(1+\varepsilon)M}{p\pi^2}\right)^{-1}\|u[0]\|, \text{ for all } t \geq 0 \tag{5.16}$$

Inequality (5.16) shows that estimate (3.8) holds with $R := G\left(1 - \frac{\gamma(1+\varepsilon)M}{p\pi^2}\right)^{-1}$ and $\omega = \delta$. The proof is complete. ◁

**Proof of Theorem 3.4:** Let an arbitrary $r \geq 0$ be given. Let also an arbitrary solution $u \in C^0(\Re_+; L^2(0,1)) \cap C^1((0,+\infty) \times [0,1])$ with $u[t] \in C^2([0,1])$ for $t>0$ of the closed-loop system (3.1), (3.2) with (3.12) be given. Define for all $t \geq 0$ and $n = 1,2,...$:

$$c_n(t) = \sqrt{2}\int_0^1 u(t,z)\sin(n\pi z)dz \tag{5.17}$$

$$v_n(t) = \sqrt{2}\int_0^1 v(t,z)\sin(n\pi z)dz \tag{5.18}$$

It follows from (3.1), (3.2), (3.12), (5.17), (5.18) and repeated integration by parts that the following differential equations hold for all $t>0$ and $n=1,2,...$:

$$\dot{c}_n(t) = -(n^2\pi^2 p + q)c_n(t) + (-1)^n n\pi^2 prc_1(t) + v_n(t) \tag{5.19}$$

Since $u \in C^0(\mathfrak{R}_+; L^2(0,1)) \cap C^1((0,+\infty) \times [0,1])$ it follows from (5.17) that $c_n \in C^0(\mathfrak{R}_+) \cap C^1((0,+\infty))$ for all $n=1,2,...$. Therefore, the following formulas hold for all $t \geq 0$:

$$c_1(t) = c_1(0)\exp\left(-\beta(\mu^2+1+r)t\right) + \int_0^t \exp\left(-\beta(\mu^2+1+r)(t-s)\right)v_1(s)ds \tag{5.20}$$

$$c_n(t) = c_n(0)\exp\left(-\beta(n^2+\mu^2)t\right) + \int_0^t \exp\left(-\beta(n^2+\mu^2)(t-s)\right)\left((-1)^n n\beta rc_1(s) + v_n(s)\right)ds,$$

$$\text{for } n = 2,3,... \tag{5.21}$$

where

$$\beta := p\pi^2, \quad \mu := \frac{1}{\pi}\sqrt{\frac{q}{p}} \tag{5.22}$$

Equations (5.20), (5.21) imply the following estimates for all $t \geq 0$:

$$|c_1(t)| \leq |c_1(0)|\exp\left(-\beta(\mu^2+1+r)t\right) + \int_0^t \exp\left(-\beta(\mu^2+1+r)(t-s)\right)|v_1(s)|ds \tag{5.23}$$

$$|c_n(t)| \leq |c_n(0)|\exp\left(-\beta(n^2+\mu^2)t\right) + n\beta r\int_0^t \exp\left(-\beta(n^2+\mu^2)(t-s)\right)|c_1(s)|ds$$

$$+ \int_0^t \exp\left(-\beta(n^2+\mu^2)(t-s)\right)|v_n(s)|ds$$

$$\text{for } n = 2,3,... \tag{5.24}$$

Using the Cauchy-Schwarz inequality and (5.23), we obtain the following estimate for all $\omega \in (0, \mu^2+1+r)$ and $t \geq 0$:

$$|c_1(t)| \leq |c_1(0)|\exp\left(-\beta(\mu^2+1+r)t\right) + \left(\frac{1}{2\beta(\mu^2+1+r-\omega)}\right)^{1/2}\left(\int_0^t \exp(-2\omega\beta(t-s))|v_1(s)|^2 ds\right)^{1/2} \tag{5.25}$$

Using the Cauchy-Schwarz inequality and (5.24), (5.25), we obtain the following estimates for all $\omega \in (0, \mu^2 + \min(4, 1+r))$, $n = 2,3,...$ and $t \geq 0$:

$$|c_n(t)| \leq |c_n(0)|\exp\left(-\beta(n^2+\mu^2)t\right)$$

$$+n\beta r|c_1(0)|\exp\left(-\beta(n^2+\mu^2)t\right)\int_0^t \exp\left(\beta(n^2-1-r)s\right)ds$$

$$+\frac{nr}{(n^2+\mu^2-\omega)}\left(\frac{1}{2\beta(1+\mu^2+r-\omega)}\right)^{1/2}\left(\int_0^t \exp(-2\omega\beta(t-s))|v_1(s)|^2 ds\right)^{1/2} \quad (5.26)$$

$$+\left(\frac{1}{2\beta(n^2+\mu^2-\omega)}\right)^{1/2}\left(\int_0^t \exp(-2\omega\beta(t-s))|v_n(s)|^2 ds\right)^{1/2}$$

Using the fact that $\exp\left(\beta(n^2-1-r)s\right) \leq \exp\left(\beta(n^2-1)s\right)$ in conjunction with (5.26), we get for all $\omega \in (0, \mu^2 + \min(4, 1+r))$, $n = 2, 3, \ldots$ and $t \geq 0$:

$$|c_n(t)| \leq \exp\left(-\beta(1+\mu^2)t\right)\left(|c_n(0)| + |c_1(0)|\frac{nr}{n^2-1}\right)$$

$$+\frac{nr}{(n^2+\mu^2-\omega)}\left(\frac{1}{2\beta(1+\mu^2+r-\omega)}\right)^{1/2}\left(\int_0^t \exp(-2\omega\beta(t-s))|v_1(s)|^2 ds\right)^{1/2} \quad (5.27)$$

$$+\left(\frac{1}{2\beta(n^2+\mu^2-\omega)}\right)^{1/2}\left(\int_0^t \exp(-2\omega\beta(t-s))|v_n(s)|^2 ds\right)^{1/2}$$

Using the fact that the inequality $(a+b)^2 \leq (1+\varepsilon^{-1})a^2 + (1+\varepsilon)b^2$ holds for all $\varepsilon > 0$, $a, b \geq 0$ in conjunction with (5.25) and the fact that $\exp\left(-\beta(\mu^2+1+r)t\right) \leq \exp\left(-\beta(\mu^2+1)t\right)$, we get for all $\varepsilon > 0$, $\omega \in (0, \mu^2 + \min(4, 1+r))$ and $t \geq 0$:

$$|c_1(t)|^2 \leq (1+\varepsilon^{-1})|c_1(0)|^2 \exp\left(-2\beta(\mu^2+1)t\right)$$

$$+\frac{1+\varepsilon}{2\beta(\mu^2+1+r-\omega)}\left(\int_0^t \exp(-2\omega\beta(t-s))|v_1(s)|^2 ds\right) \quad (5.28)$$

Using the fact that the inequality $(a+b)^2 \leq (1+\varepsilon^{-1})a^2 + (1+\varepsilon)b^2$ holds for all $\varepsilon > 0$, $a, b \geq 0$ in conjunction with (5.27), the fact that $\left(|c_n(0)| + |c_1(0)|\frac{nr}{n^2-1}\right)^2 \leq 2|c_n(0)|^2 + 2|c_1(0)|^2\frac{n^2r^2}{(n^2-1)^2}$ as well as the fact that

$$A^2 \leq \frac{n^2r^2}{(n^2+\mu^2-\omega)^2}\left(\frac{1+k_n^{-1}}{2\beta(1+\mu^2+r-\omega)}\right)\left(\int_0^t \exp(-2\omega\beta(t-s))|v_1(s)|^2 ds\right)$$

$$+\left(\frac{1+k_n}{2\beta(n^2+\mu^2-\omega)}\right)\left(\int_0^t \exp(-2\omega\beta(t-s))|v_n(s)|^2 ds\right)$$

where

$$A = \frac{nr}{(n^2+\mu^2-\omega)}\left(\frac{1}{2\beta(1+\mu^2+r-\omega)}\right)^{1/2}\left(\int_0^t \exp(-2\omega\beta(t-s))|v_1(s)|^2 ds\right)^{1/2}$$

$$+\left(\frac{1}{2\beta(n^2+\mu^2-\omega)}\right)^{1/2}\left(\int_0^t \exp(-2\omega\beta(t-s))|v_n(s)|^2 ds\right)^{1/2}$$

$$k_n = \frac{(n^2+\mu^2-\omega)-\lambda(4+\mu^2-\omega)}{\lambda(4+\mu^2-\omega)}$$

and $\lambda \in (0,1)$, gives the following estimates for all $\varepsilon > 0$, $\lambda \in (0,1)$, $\omega \in (0, \mu^2 + \min(4, 1+r))$, $n = 2,3,...$ and $t \geq 0$:

$$|c_n(t)|^2 \leq 2(1+\varepsilon^{-1})\exp(-2\beta(1+\mu^2)t)\left(|c_n(0)|^2 + |c_1(0)|^2 \frac{n^2 r^2}{(n^2-1)^2}\right)$$

$$+ \frac{(1+\varepsilon)n^2 r^2 \left(\int_0^t \exp(-2\omega\beta(t-s))|v_1(s)|^2 ds\right)}{2\beta(1+\mu^2+r-\omega)(n^2+\mu^2-\omega)\left(n^2+\mu^2-\omega-\lambda(4+\mu^2-\omega)\right)} \quad (5.29)$$

$$+ \frac{1+\varepsilon}{2\beta\lambda(4+\mu^2-\omega)}\left(\int_0^t \exp(-2\omega\beta(t-s))|v_n(s)|^2 ds\right)$$

Since the set of functions $\phi_n(x) = \sqrt{2}\sin(n\pi x)$, $n = 1,2,3,...$ is an orthonormal basis of $L^2(0,1)$, Parseval's identity in conjunction with definitions (5.17), (5.18) gives for all $t \geq 0$:

$$\|u[t]\|^2 = \sum_{n=1}^\infty c_n^2(t) \quad (5.30)$$

$$\|v[t]\|^2 = \sum_{n=1}^\infty v_n^2(t) \quad (5.31)$$

Combining (5.30), (5.29) and (5.28), we get for all $\varepsilon > 0$, $\lambda \in (0,1)$, $\omega \in (0, \mu^2 + \min(4, 1+r))$ and $t \geq 0$:

$$\|u[t]\|^2 \leq 2(1+\varepsilon^{-1})\exp(-2\beta(1+\mu^2)t)\left(\|u[0]\|^2 + |c_1(0)|^2 r^2 \sum_{n=2}^\infty \frac{n^2}{(n^2-1)^2}\right)$$

$$+\left(1 + r^2 \sum_{n=2}^\infty \frac{n^2}{(n^2+\mu^2-\omega)(n^2+\mu^2-\omega-\lambda(4+\mu^2-\omega))}\right)\frac{(1+\varepsilon)\int_0^t \exp(-2\omega\beta(t-s))|v_1(s)|^2 ds}{2\beta(1+\mu^2+r-\omega)} \quad (5.32)$$

$$+ \frac{1+\varepsilon}{2\beta\lambda(4+\mu^2-\omega)}\sum_{n=2}^\infty \int_0^t \exp(-2\omega\beta(t-s))|v_n(s)|^2 ds$$

Using the fact that $|c_1(0)| \leq \|u[0]\|$, the fact that $\frac{n^2}{n^2-1} \leq \frac{4}{3}$ for all $n = 2,3,...$, (5.31) and (5.32), we obtain for all $\varepsilon > 0$, $\lambda \in (0,1)$, $\omega \in (0, \mu^2 + \min(4, 1+r))$ and $t \geq 0$:

$$\|u[t]\|^2 \le 2(1+\varepsilon^{-1})\exp(-2\beta(1+\mu^2)t)\left(1+\frac{4r^2}{3}\sum_{n=2}^{\infty}\frac{1}{n^2-1}\right)\|u[0]\|^2$$
$$+\frac{(1+\varepsilon)K}{2\beta(4+\mu^2-\omega)}\int_0^t \exp(-2\omega\beta(t-s))\|v[s]\|^2\,ds \tag{5.33}$$

where

$$K:=\lambda^{-1}\max\left(\frac{\lambda(4+\mu^2-\omega)}{1+\mu^2+r-\omega}\left(1+r^2\sum_{n=2}^{\infty}\frac{n^2}{(n^2+\mu^2-\omega)(n^2+\mu^2-\omega-\lambda(4+\mu^2-\omega))}\right),1\right) \tag{5.34}$$

Using the fact that $\frac{1}{n^2-1}\le\frac{4}{3n^2}$ for all $n=2,3,\ldots$, the fact that $\sum_{n=1}^{\infty}\frac{1}{n^2}=\frac{\pi^2}{6}$, we obtain from (5.33) for all $\varepsilon>0$, $\lambda\in(0,1)$, $\omega\in(0,\mu^2+\min(4,1+r))$, $\sigma\in[0,\omega)$ and $t\ge 0$:

$$\|u[t]\|^2 \le 2(1+\varepsilon^{-1})\exp(-2\beta(1+\mu^2)t)\left(1+\frac{8(\pi^2-6)r^2}{27}\right)\|u[0]\|^2$$
$$+\frac{(1+\varepsilon)K}{4\beta^2(\omega-\sigma)(4+\mu^2-\omega)}\sup_{0\le s\le t}\left(\exp(-2\sigma\beta(t-s))\|v[s]\|^2\right) \tag{5.35}$$

Define $\tilde{\omega}:=\omega-\mu^2$ and notice that definitions (5.34), (3.15) imply that $K=L(r,\tilde{\omega},\lambda)$. Inequality (3.13) with $\tilde{\omega}:=\omega-\mu^2$ in place of $\omega$ is a direct consequence of (5.35), definitions (3.14), (5.22) and the fact that $\sqrt{a+b}\le\sqrt{a}+\sqrt{b}$ for all $a,b\ge 0$. The proof is complete. ◁

**Proof of Theorem 4.1:** Consider the transformation

$$w(t,x)=u(t,x)-\sum_{i=1}^{n}\varphi_i(x)\xi_i(t),\text{ for }t\ge 0,\ x\in[0,1] \tag{5.36}$$

Using (4.1), (4.6) and (5.36), we conclude that the following PDE holds:

$$w_t(t,x)=pw_{xx}(t,x)-qw(t,x)+(f(u[t]))(x)-\sum_{i=1}^{n}\varphi_i(x)K_i(u[t])+\sum_{i=1}^{n}\left(p\varphi_i''(x)-(q-\omega_i^2)\varphi_i(x)\right)\xi_i(t),$$
$$\text{for }t>0,\ x\in(0,1) \tag{5.37}$$

Using (3.2), (4.5), (5.36) and the fact that $\varphi_i(0)=0$ for $i=1,\ldots,n$, we also obtain the following equations:

$$w(t,0)=0,\text{ for }t\ge 0 \tag{5.38}$$

$$w(t,1)=\int_0^1 k(x)w(t,x)dx,\text{ for }t\ge 0 \tag{5.39}$$

It follows from (5.36), (5.37), (5.38), (5.39) that for every solution $u\in C^0(\Re_+;L^2(0,1))\cap C^1((0,+\infty)\times[0,1])$, $\xi_i\in C^0(\Re_+)\cap C^1((0,+\infty))$, $i=1,\ldots,n$ with $u[t]\in C^2([0,1])$ for $t>0$ of the closed-loop system (4.1), (3.2) with (4.5), (4.6), the function $w$ is of class $C^0(\Re_+;L^2(0,1))\cap C^1((0,+\infty)\times[0,1])$ with $w[t]\in C^2([0,1])$ for $t>0$, which satisfies (5.38), (5.39) as well as the PDE

$$w_t(t,x)=pw_{xx}(t,x)-qw(t,x)+v(t,x),\text{ for }t>0,\ x\in(0,1) \tag{5.40}$$

with

$$v(t,x) := (f(u[t]))(x) - \sum_{i=1}^{n} \varphi_i(x) K_i(u[t]) + \sum_{i=1}^{n} \left( p\varphi_i''(x) - (q - \omega_i^2)\varphi_i(x) \right) \xi_i(t), \text{ for } t \geq 0, \ x \in [0,1]. \quad (5.41)$$

Since the kernel $k \in C^0([0,1])$ achieves Input-to-State Stabilization with gain $\frac{\gamma}{p\pi^2} \geq 0$ of system (3.1), (3.2), there exist constants $G, \sigma > 0$ so that every solution $w \in C^0(\Re_+; L^2(0,1)) \cap C^1((0,+\infty) \times [0,1])$ with $w[t] \in C^2([0,1])$ for $t > 0$ of the system (5.40), (5.38), (5.39) satisfies the estimate:

$$\|w[t]\| \leq G \exp(-\sigma t) \|w[0]\| + \frac{\gamma}{p\pi^2} \sup_{0 \leq s \leq t} \left( \|v([s])\| \right), \text{ for all } t \geq 0 \quad (5.42)$$

Lemma 7.1 in [15] in conjunction with inequality (5.42) implies that for every $\varepsilon > 0$ there exists a constant $\delta \in (0, \sigma)$ with the following property: every solution $w \in C^0(\Re_+; L^2(0,1)) \cap C^1((0,+\infty) \times [0,1])$ with $w[t] \in C^2([0,1])$ for $t > 0$ of the closed-loop system (5.40), (5.38), (5.39) satisfies the estimate:

$$\|w[t]\| \leq G \exp(-\delta t) \|w[0]\| + \frac{\gamma(1+\varepsilon)}{p\pi^2} \sup_{0 \leq s \leq t} \left( \|v([s])\| \exp(-\delta(t-s)) \right), \text{ for all } t \geq 0 \quad (5.43)$$

Inequality (4.4) allows us to pick sufficiently small $\varepsilon > 0$ so that (5.43) holds as well as the following inequality holds:

$$L := \frac{\gamma(1+\varepsilon)}{p\pi^2} \left( M + \sum_{i=1}^{n} \frac{P_i \|p\varphi_i'' - (q - \omega_i^2)\varphi_i\|}{\omega_i^2 - \varepsilon} \right) + \sum_{i=1}^{n} \frac{P_i \|\varphi_i\|}{\omega_i^2 - \varepsilon} < 1 \quad (5.44)$$

Without loss of generality, we may assume that $\delta \leq \varepsilon < \omega_i^2$ for $i = 1,\ldots,n$. Using (4.3), (5.41), the triangle inequality, we conclude from (5.43) that the following estimate holds for all $t \geq 0$:

$$\|w[t]\| \leq G \exp(-\delta t) \|w[0]\| + \frac{\gamma(1+\varepsilon)}{p\pi^2} M \max_{0 \leq s \leq t} \left( \exp(-\delta(t-s)) \|u[s]\| \right)$$
$$+ \frac{\gamma(1+\varepsilon)}{p\pi^2} \sum_{i=1}^{n} \max_{0 \leq s \leq t} \left( \exp(-\delta(t-s)) \|p\varphi_i'' - (q - \omega_i^2)\varphi_i\| |\xi_i(s)| \right) \quad (5.45)$$

Using (4.6) and the variations of constants formula, we obtain the following estimates for $i = 1,\ldots,n$ and $t \geq 0$:

$$|\xi_i(t)| \leq \exp(-\omega_i^2 t) |\xi_i(0)| + \frac{1}{\omega_i^2 - \delta} \sup_{0 \leq s \leq t} \left( |K_i(u[s])| \exp(-\delta(t-s)) \right) \quad (5.46)$$

Combining (4.2) and (5.46), we obtain the following estimates for $i = 1,\ldots,n$ and $t \geq 0$:

$$|\xi_i(t)| \leq \exp(-\omega_i^2 t) |\xi_i(0)| + \frac{P_i}{\omega_i^2 - \delta} \max_{0 \leq s \leq t} \left( \|u[s]\| \exp(-\delta(t-s)) \right) \quad (5.47)$$

Next, we use definition (5.36) and the triangle inequality, which implies the following estimates for $t \geq 0$:

$$\|u[t]\| \leq \|w[t]\| + \sum_{i=1}^{n} \|\varphi_i\| |\xi_i(t)| \quad (5.48)$$

Combining (5.44), (5.45), (5.47), (5.48) and using the fact that $\delta \leq \varepsilon < \omega_i^2$ for $i = 1,\ldots,n$, we obtain for $t \geq 0$:

$$\|w[t]\|\exp(\delta t) \le G\|w[0]\| + \sum_{i=1}^{n}\left(\frac{\gamma(1+\varepsilon)}{p\pi^2}\|p\varphi_i'' - (q-\omega_i^2)\varphi_i\| + \|\varphi_i\|\right)|\xi_i(0)|$$
$$+ L\max_{0\le s\le t}\left(\exp(\delta s)\|u[s]\|\right) \quad (5.49)$$

It follows from (5.44) and (5.49) that there exists a constant $\Theta > 0$ such that the following inequality holds for all $t \ge 0$:

$$\|u[t]\|\exp(\delta t) \le \Theta\left(\|w[0]\| + \sum_{i=1}^{n}|\xi_i(0)|\right) \quad (5.50)$$

Inequality (4.7) for appropriate constant $R > 0$ is a direct consequence of inequality (5.50) and the fact that $\|w[0]\| \le \|u[0]\| + \sum_{i=1}^{n}\|k_i\||\xi_i(0)|$ (a consequence of definition (5.36) and the triangle inequality). The proof is complete.   ◁

## 6. Back to the Motivating Example

Here we turn again our attention to the stabilization problem for system (2.1), (2.2) with

$$p = 1, \quad (f(u))(x) = \|u\|\sum_{i=1}^{n}A_i\varphi_i(x) \text{ for } x \in [0,1] \quad (6.1)$$

where $A_i \in \Re$, $\omega_i > 0$, $i = 1,...,n$ are constants and the functions $\varphi_i(x)$, $i = 1,...,n$ are given by:
- $\varphi_i(x) = x$, when $q = \omega_i^2$,
- $\varphi_i(x) = \sin\left(x\sqrt{\omega_i^2 - q}\right)$, when $q < \omega_i^2$,
- $\varphi_i(x) = \sinh\left(x\sqrt{q - \omega_i^2}\right)$, when $q > \omega_i^2$,

The basic boundary feedback kernel $k \in C^0([0,1])$ that achieves Input-to-State Stabilization with gain $\pi^{-2}\gamma \ge 0$ of system (3.1), (3.2) may be selected as follows:
- when $q \ge 0$ then we may select the single-mode boundary feedback $k(x) = -\pi r\sin(\pi x)$ with $r \ge 0$. The fact that this kernel achieves Input-to-State Stabilization of system (3.1), (3.2) is a consequence of Theorem 3.4.
- when $q < 0$ then we may select a boundary feedback designed by means of backstepping (see [17,24]) and take $k(x) = -(r-q)x\dfrac{I_1\left(\sqrt{(r-q)(1-x^2)}\right)}{\sqrt{(r-q)(1-x^2)}}$ with $r \ge 0$, where $I_1$ is the modified Bessel function of the first kind of order one. The fact that this kernel achieves Input-to-State Stabilization of system (3.1), (3.2) follows from an analysis similar to that given in [12,15].

Selecting $K_i(u) = A_i\|u\|$, we conclude that inequalities (4.2), (4.3) hold with $M = 0$, $P_i = |A_i|$, $i = 1,...,n$. Therefore, Theorem 4.1 implies that $0 \in L^2(0,1) \times \Re^n$ is globally exponentially stable for the closed-loop system (2.1), (2.2) with

$$U(t) = \sum_{i=1}^{n}\varphi_i(1)\xi_i(t) + \int_0^1 k(x)u(t,x)dx - \sum_{i=1}^{n}\xi_i(t)\int_0^1 k(x)\varphi_i(x)dx \quad (6.2)$$

$$\dot{\xi}_i(t) = -\omega_i^2\xi_i(t) + A_i\|u[t]\|, \quad i = 1,...,n \quad (6.3)$$

provided that the small-gain condition

$$\sum_{i=1}^{n} \omega_i^{-2} |A_i| \|\varphi_i\| < 1 \tag{6.4}$$

holds.

It is important to compare condition (6.4) with the condition arising for the case of static boundary feedback. Proposition 3.2 indicates that that $0 \in L^2(0,1)$ is globally exponentially stable for the closed-loop system (2.1), (2.2) with

$$U(t) = \int_0^1 k(x)u(t,x)dx \tag{6.5}$$

provided that the small-gain condition

$$\frac{\gamma}{\pi^2} \left\| \sum_{i=1}^{n} A_i \varphi_i(x) \right\| < 1 \tag{6.6}$$

holds. Notice that (6.4) is *independent* of the gain $\frac{\gamma}{p\pi^2} \geq 0$ achieved by the basic boundary feedback kernel $k \in C^0([0,1])$, while (6.6) depends heavily on this gain. Moreover, there are cases where (6.4) is much less demanding than (6.6). For example, when $q = 0$, $n = 1$, i.e., when (2.3) holds, condition (6.4) is equivalent to the condition

$$|A| < \frac{2\omega^{5/2}}{\sqrt{2\omega - \sin(2\omega)}} \tag{6.7}$$

while condition (6.6) in conjunction with Theorem 3.3 implies the condition

$$|A| < \frac{2\pi^3 \omega^{1/2}}{\sqrt{\pi^2 - 6}\sqrt{2\omega - \sin(2\omega)}} \tag{6.8}$$

It is clear that when $\omega^2 > \frac{\pi^3}{\sqrt{\pi^2 - 6}}$, condition (6.7) is less demanding than condition (6.8).

Indeed, for system (2.1), (2.2), (2.3) with $\omega = 20$ and $A = 500$, condition (6.8) does not hold. Consequently, Theorem 3.3 guarantees that condition (6.6) will not hold for any possible selection of the basic boundary feedback kernel $k \in C^0([0,1])$. On the other hand, condition (6.7) holds and consequently for every $r \geq 0$, the equilibrium point $0 \in L^2(0,1) \times \Re$ is globally exponentially stable for the closed-loop system (2.1), (2.2), (2.3) with $\omega = 20$, $A = 500$ and

$$U(t) = -\pi r \int_0^1 \sin(\pi x) u(t,x)dx + \xi(t)\left( \sin(\omega) + \pi r \frac{\sin(\omega - \pi)}{2(\omega - \pi)} - \pi r \frac{\sin(\omega + \pi)}{2(\omega + \pi)} \right) \tag{6.9}$$

$$\dot{\xi}(t) = -\omega^2 \xi(t) + A\|u[t]\| \tag{6.10}$$

Figure 3 and Figure 4 show the evolution of the state norm $\|u[t]\|$ for the same initial condition and different values of $r \geq 0$. An increase of $r \geq 0$ causes a faster convergence but with a higher overshoot.

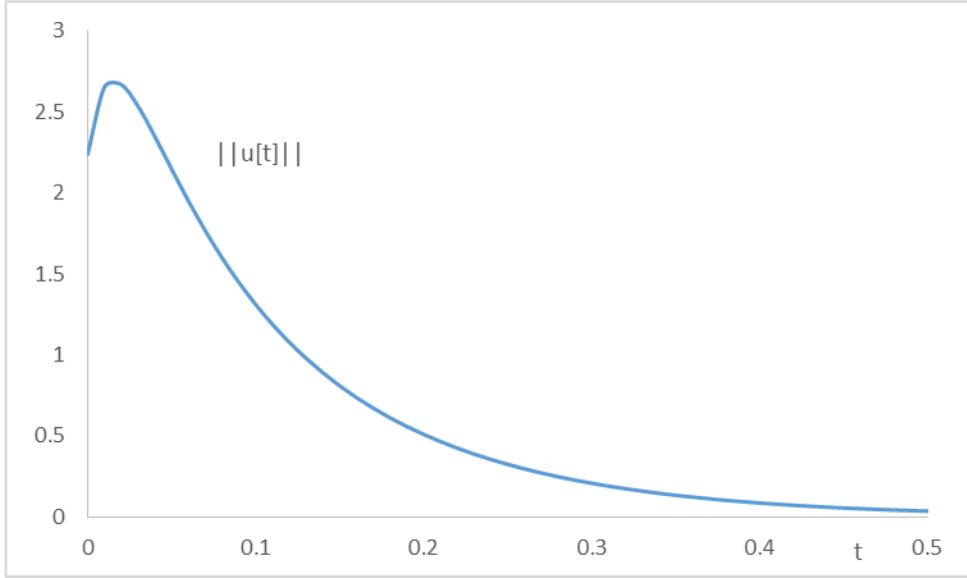

**Fig. 3:** The evolution of the state norm $\|u[t]\|$ for the closed-loop system (2.1), (2.2), (2.3), (6.9), (6.10) with $r=0$, $\omega=20$, $A=500$. The initial condition is $\xi(0)=0$,
$$u_0(x) = \sqrt{2}\sin(\pi x) + 2\sqrt{2}\sin(2\pi x) + \frac{\sqrt{2}}{10}\sin(3\pi x) + \frac{\sqrt{2}}{100}\sum_{n=37}^{42}\sin(n\pi x),\ x\in[0,1].$$

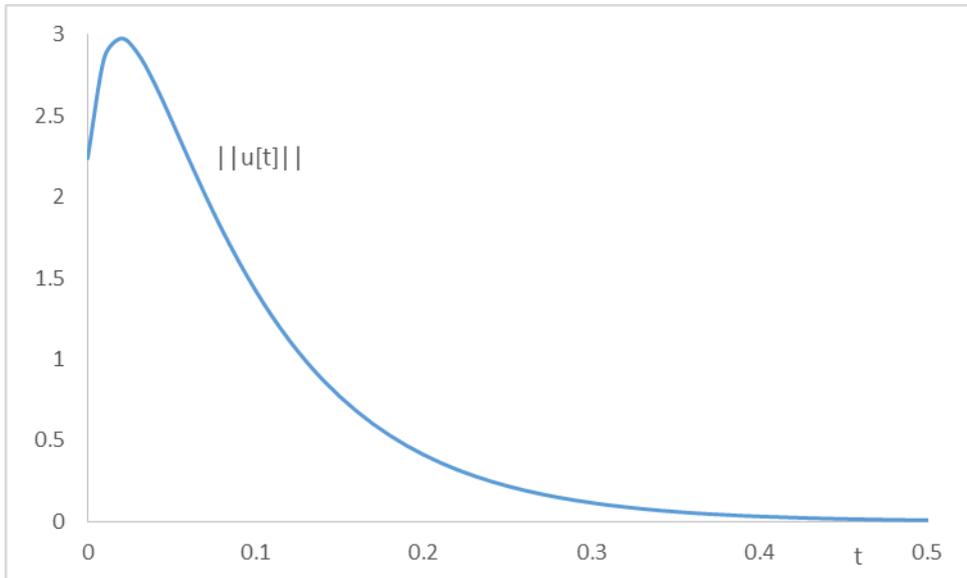

**Fig. 4:** The evolution of the state norm $\|u[t]\|$ for the closed-loop system (2.1), (2.2), (2.3), (6.9), (6.10) with $r=0.9$, $\omega=20$, $A=500$. The initial condition is $\xi(0)=0$,
$$u_0(x) = \sqrt{2}\sin(\pi x) + 2\sqrt{2}\sin(2\pi x) + \frac{\sqrt{2}}{10}\sin(3\pi x) + \frac{\sqrt{2}}{100}\sum_{n=37}^{42}\sin(n\pi x),\ x\in[0,1].$$

# 7. Concluding Remarks

The present paper introduces a small-gain methodology for the design of boundary feedback stabilizers in 1-D, semilinear, parabolic PDEs with nonlocal terms. The stabilization results are global and the nonlinearities are assumed to satisfy a linear growth condition with restricted linear growth coefficient. Two different boundary feedback stabilizers are provided: a linear static boundary feedback stabilizer, which can handle uncertain nonlinearities and a nonlinear dynamic boundary feedback stabilizer. However, the paper also provides additional results which reveal fundamental limitations for feedback design in the parabolic case: the fact that gain assignment is not possible by means of boundary feedback.

The focus of the paper is on the specific case of Dirichlet actuation at one end of the domain and Dirichlet boundary condition at the other end. Future research will address the extension to other cases (Neumann actuation, Robin or Neumann boundary conditions at the non-actuated end of the domain).